\documentclass[a4paper,10pt]{amsart}
\usepackage{amsmath,amscd,amssymb,epsfig}

\newtheorem{theo}{\bf Theorem}

\newtheorem{coro}[theo]{\bf Corollary}

\newenvironment{dem}[1][\bf Proof]
  {\begin{trivlist}\item{\em #1:}\ }
  {\hfill$\square$\par\end{trivlist}}
\newenvironment{dem...}[1][\bf Proof]
  {\begin{trivlist}\item{\em #1:}\ }
  {\par\end{trivlist}}

\newenvironment{Rq}[1][\bf Remark]
  {\begin{trivlist}\item{\em #1:}\ }
  {\par\end{trivlist}}

\def\NN{\mathbb{N}}  
\def\ZZ{\mathbb{Z}}  
\def\QQ{\mathbb{Q}}
\def\L{\Lambda}
\def\A{{\mathcal{A}}}
\def\B{{\mathcal{B}}}
\def\go{\longrightarrow}
\def\S{\mathfrak{S}}
\def\co{\colon\thinspace}
\newcommand{\color}[6]{}

\begin{document}
\title{Non injectivity of the ``hair'' map}
\author{Bertrand Patureau-Mirand}
\address{LMAM, universit\'e de Bretagne-Sud, universit\'e europ\'eenne de
  Bretagne, BP 573, 56017 Vannes, France }
\email{bertrand.patureau@univ-ubs.fr}
\date{\today}

\begin{abstract}
  Kricker constructed a knot invariant $Z^{rat}$ valued in a space of
  Feynman diagrams with beads.  When composed with the so called
  ``hair'' map $H$, it gives the Kontsevich integral of the knot. We
  introduce a new grading on diagrams with beads and use it to show
  that a non trivial element constructed from Vogel's zero divisor in
  the algebra $\L$ is in the kernel of $H$. This shows that $H$ is not
  injective.
\end{abstract}

\maketitle

\section*{Introduction}
The Kontsevich integral $Z$ is a universal rational finite type
invariant for knots (see the Bar-Natan survey \cite{BN}). For a knot
$K$, $Z(K)$ lives in the space of Chinese diagrams isomorphic to
$\hat\B(*)$ (see Section 1.1). Rozansky conjectured (\cite{Roz}) and
Kricker proved (\cite{Kr}) that $Z$ can be organized into a series of
"lines" called $Z^{rat}$. They can be represented by finite
$\QQ$--linear combinations of diagrams whose edges are labelled, in an
appropriate way, with rational functions. In \cite{GaKr}, Garoufalidis
and Kricker directly proved that the map $Z^{rat}$ with values in a
space of diagrams with beads is an isotopy invariant and that $Z$
factors through $Z^{rat}$. For a knot $K$ with trivial Alexander
polynomial, $Z(K)=H\circ Z^{rat}(K)$ where $H$ is the hair map (see
Section 1.3).  Rozansky, Garoufalidis and Kricker conjectured (see
\cite[Conjecture 3.18]{Oth}) that $H$ could be injective. Theorem
\ref{T:main} gives a counterexample to this conjecture.

\section{The hair map}
\subsection{Classical diagrams}
Let $X$ be a finite set.  A $X$--diagram is an isomorphism class of finite
uni-tri-valent graphs $K$ with the following data:
\begin{itemize}
\item At each trivalent vertex $x$ of $K$, we have a cyclic ordering
  on the three oriented edges starting from $x$.
\item A bijection between the set of univalent vertices of $K$ and the set
  $X$.
\end{itemize}
We define $\A(X)$ to be the quotient of the $\QQ$--vector space generated by
$X$--diagrams by the relations:
\begin{enumerate}
\item The (AS) relations for
  ``antisymmetry'': $$\begin{array}{cccc}\put(-10,-10)
    {\epsfbox{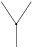}}&+& \put(-10,-10)
    {\epsfbox{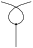}}&=0\end{array}$$
\item The (IHX) relations for three diagrams which differ only in a
neighborhood of an edge:
  $$\begin{array}{ccccc}\put(-8,-5) {\epsfbox{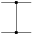}}&=&
    \put(-8,-5) {\epsfbox{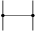}} &
    -&\put(-8,-5) {\epsfbox{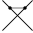}}\end{array}$$
\end{enumerate}
These spaces are graded.  The degree of an $X$--diagram is given by half the
total number of vertices.\\ 
Let $[n]=\{1,2,\ldots,n\}$ and define $F_n$ to be the subspace of $\A([n])$
generated by connected diagrams with at least one trivalent vertex. The
permutation group $\S(X)$ acts on $\A(X)$. Let $B(*)$ be the coinvariant space
for this action:
$$B(*)=\bigoplus_{n\in\NN}\A([n])\otimes_{\S_n}\QQ$$
and let $\hat B(*)$ be the completion of $B(*)$ for the grading.\\
Finally let $\L$ be Vogel's algebra generated by totally antisymmetric
elements of $F_3$ (for the action of $\S_3$).\\
We recall (see \cite{Vo}) that $\L$ acts on the modules $F_n$ and that for this
action, $F_0$ and $F_2$ are free $\L$--modules of rank one. Furthermore, the
following elements are in $\L$:
$$\begin{array}{cccccc}t=&\put(-8,-5) {\epsfbox{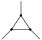}}&=\frac12&
  \put(-8,-5) {\epsfbox{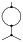}} & \qquad x_n=&\put(-8,-10) 
  {\epsfbox{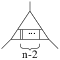}}
\end{array}$$

\begin{theo}(Vogel \cite[Section 8 and Proposition 8.5]{Vo}). The
  element $t$ is a divisor of zero in $\Lambda$.
\end{theo}

\begin{coro}\label{C}
There exists an element $r\in\L\setminus\{0\}$ such that $t.r=0$. So one has
$$\put(-90,-15){\input{t-killedF0.pstex_t}}\quad\hbox{but}
\quad\put(20,-15){\input{t-killedF3.pstex_t}}$$
\end{coro}
\begin{dem}
  $F_0$ is a free $\L$--module of rank one generated by the diagram
  $\Theta$ and the previous diagram of $F_0$ is $r.\Theta\neq0$. The
  diagram of $F_3$ of the corollary is the product
  $r.\put(2,-7){\epsfbox{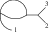}}\hspace{25pt}=2tr=0\in\L$.
\end{dem}

\begin{Rq} 
  Vogel shows that $r$ can be chosen with degree fifteen in $\L$ (the
  degree in $\L$ is the degree in $F_3$ minus two), and in the algebra
  generated by the $x_n$.  This element is killed by all the weight
  systems coming from Lie algebras (but $r$ is not killed by the Lie
  superalgebras $\mathfrak{D}_{2,1,\alpha}$).
\end{Rq}

\subsection{Diagrams with beads}
Diagrams with beads have been introduced by Kricker and Garoufalidis (see
\cite{Kr}, \cite{GaKr}). A presentation of $\B$ which uses the first
cohomology classes of diagrams is already present in \cite{Roz}. Vogel
explained me this point of view for diagrams with beads.\\

Let $G$ be the multiplicative group $\{b^n,\,n\in\ZZ\}\simeq (\ZZ,+)$ and
consider its group algebra $R=\QQ G =\QQ[b,b^{-1}]$. Let $a\mapsto \overline
a$ be the involution of the $\QQ$--algebra $R$ that maps $b$ to $b^{-1}$.
\\
A diagram with beads in $R$ is an $\emptyset$--diagram with the following
supplementary data: The beads form a map $f\co E\go R$ from the set of oriented
edges of $K$ such that if $-e$ denotes the same edge than $e$ with opposite
orientation, one has $f(-e)=\overline{f(e)}$.
\\
We will represent the beads by some arrows on the edges with label in $R$. The
value of the bead $f$ on $e$ is given by the product of these labels and we
will not represent the beads with value $1$. So with graphical notations, we
have:
$$\input{involution.pstex_t}\quad\hbox{and}\quad\input{Multilin1.pstex_t}$$
The loop degree of a diagram with beads is the first Betti number of the
underlying graph.
\\
Let $\A^R(\emptyset)$ be the quotient of the $\QQ$--vector space generated
by diagrams with beads in $R$ by the following relations:
\begin{enumerate}
\item (AS) 
\item The (IHX) relations should only be considered near an edge with bead $1$.
\item PUSH:$$\put(-50,-10){\epsfbox{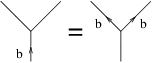}}$$
\item Multilinearity: $$\input{Multilin2.pstex_t}$$
\end{enumerate}
$\A^R(\emptyset)$ is graded by the loop degree:
$$\A^R(\emptyset)=\bigoplus_{n\in\NN}\A^R_n(\emptyset)$$
We will prefer another presentation of $\A^R(\emptyset)$:
\begin{itemize}
\item Remark that it is enough to consider diagrams with beads in $G$
  and the multilinear relation can be viewed as a notation.
\item Next remark that for a diagram with beads in $G$, the map $f$
  define a $1$--cochain $\tilde f$ with values in $\ZZ\simeq G$ on the
  underlying simplicial set of $K$. The elements $\tilde f$ are in
  fact $1$--cocycles because of the condition $f(-e)=\overline{f(e)}$
  which implies $\tilde f(-e)=-\tilde f(e)$.
\item The ``PUSH'' relation at a vertex $v$ implies that $\tilde f$
  is only given up to the coboundary of the $0$--cochain with value $1$
  on $v$ and $0$ on the other vertices. Hence $\A^R(\emptyset)$ is
  also the $\QQ$--vector space generated by the pairs ($3$--valent graph
  $D$, $x\in H^1(D,\ZZ)$) quotiented by the relations (AS) and
  (IHX). With these notations one can describe the (IHX) relations in
  the following way:
  \\
  Let $K_I$, $K_H$ and $K_X$ be three graphs which appear in a (IHX) relation
  on an edge $e$. Let $K_\bullet$ be the graph obtained by collapsing the edge
  $e$.  The maps $p_?\co K_?\go K_\bullet$ induce three cohomology isomorphisms.
  \\
  If $x\in H^1(K_\bullet,\ZZ)$ then the (IHX) relation at $e$ says that
  $$(K_I,p_I^*x)=(K_H,p_H^*x)-(K_X,p_X^*x)$$ holds in $\A^R(\emptyset)$.
\end{itemize}

\subsection{The hair map}
The hair map $H\co\A^R(\emptyset)\go\hat\B(*)$ replaces beads by legs (or hair):
Just replace a bead $b^n$ by the exponential of $n$ times a leg.
$$\input{hairmap.pstex_t}$$
$H$ is well defined (see \cite{GaKr}).

\section{Grading on diagrams with beads}
Remark that for a $3$--valent graph $K$, $H^1(K,\ZZ)$ is a free
$\ZZ$--module. The beads $x\in H^1(K,\ZZ)$ which occur in an (AS) or
(IHX) relation are the same up to isomorphisms. We will call $p\in\NN$
the bead degree of $(K,x)$ if $x$ is $p$ times an indivisible element
of $H^1(K,\ZZ)$.
\begin{theo}
The bead degree is well defined in $\A^R_n(\emptyset)$. Thus we have a grading
$$\A^R_n(\emptyset)=\bigoplus_{p\in\NN} \A^R_{n,p}(\emptyset)$$ where
$\A^R_{n,p}(\emptyset)$ is the subspace of $\A^R_n(\emptyset)$ generated by
diagrams with bead degree $p$.\\
Furthermore, $\A^R_{n,0}(\emptyset)\simeq \A_n(\emptyset)$ and for $p>0$,
$\A^R_{n,p}(\emptyset)\simeq\A^R_{n,1}(\emptyset)$.
\end{theo}
\begin{dem}
  The second presentation we have given for $\A^R_n(\emptyset)$
  implies that this degree is well defined.  Indeed, the elements in a
  IHX relation have the same degree because the set of indivisible
  elements of the cohomology is preserved by isomorphisms.
  \\
  Now, the map $\psi\co R\go \QQ$ that sends $b$ to $1$ induces the
  isomorphism $\A^R_{n,0}(\emptyset)\simeq \A_n(\emptyset)$ and the
  group morphism $\phi_p\co G\go G$ that sends $b$ to $b^p$ (or the
  multiplication by $p$ in $H^1(.,\ZZ)$) induces the isomorphism
  $\A^R_{n,1}(\emptyset)\simeq\A^R_{n,p}(\emptyset)$. These maps are
  isomorphisms because they have obvious inverses.  
  %
\end{dem}

\section{A non trivial element in the kernel of $H$}

\begin{theo}\label{T:main}
  This non trivial element of $\A^R(\emptyset)$ is in the kernel of $H$:
  $$\input{hairkilled.pstex_t}$$
  Thus $H$ is not injective.
\end{theo}
\begin{dem}
  This element is not zero because its bead degree zero part is the
  opposite of the element $r.\Theta$ of Corollary \ref{C}.  Then, one has
  $$\input{hairkilled2.pstex_t}$$ 
  but all these diagrams are zero in $B(*)$ because they contain, as a
  sub-diagram, the element of $F_3$ of Corollary \ref{C}.
\end{dem}

\begin{Rq}
  The element of theorem \ref{T:main} has a loop degree seventeen.
  \\
  The hair map is obviously injective on the space of diagrams with bead
  degree zero. I don't know if the same is true in other degrees.
  \\
\end{Rq}

\linespread{1}

\end{document}